\documentclass[a4paper,10pt]{article}

\usepackage[left=2.5cm, right=2.5cm,
top=2.0cm, bottom=3.0cm]{geometry}

\usepackage[belowskip=-3pt,aboveskip=3pt,font=footnotesize]{caption}  
\usepackage[utf8]{inputenc}
 
\usepackage[compact]{titlesec}
\titlespacing{\section}{0pt}{*1.5}{*0.5}
\titlespacing{\subsection}{0pt}{*0.75}{*0.25}
\usepackage{graphicx}
\usepackage{multicol}
\usepackage{wrapfig}
\usepackage{amsmath}
\usepackage{amsthm}
\usepackage{amssymb}
\usepackage{bbold}
\usepackage{amsfonts}
\usepackage{multirow}
\usepackage{color}

\usepackage{tikz}
\usepackage{pgfplots}

\graphicspath{{./figures/}}

\newcommand{\R}{{\mathbb{R}}}
\newcommand{\N}{{\mathbb{N}}}

\newcommand{\diff}{{\,\mbox{d}}}

\newcommand{\diffEps}{{\tilde\epsilon}}
\newcommand{\sharpEps}{{\epsilon}}
\newcommand{\diffEnergy}{{E_\diffEps}}
\newcommand{\sharpEnergy}{{F_\sharpEps}}

\newcommand{\penaltyConst}{{\sigma}}
\newcommand{\x}{x}
\newcommand{\y}{y}
\newcommand{\Y}{{\textbf{Y}}}
\newcommand{\upper}{{C}}

\newcommand{\numIk}{{I_k}}
\newcommand{\stripe}{\mathcal{S}}
\newcommand{\area}{\mathcal{R}}
\newcommand{\spike}{{s}}
\newcommand{\trunk}{{t}}
\newcommand{\facet}{{f}}

\newcommand{\auto}{{\mathcal{A}}}
\newcommand{\per}{{\mathcal{P}}}
\newcommand{\branch}{{\mathcal{B}}}
\newcommand{\comp}{{\mathcal{C}}}
\newcommand{\word}{{\omega}}

\newcommand{\dxdy}{\diff x \diff y }

\newcommand{\elastE}{{F_{\text{elast}}}}
\newcommand{\elastEK}{{F_{\text{elast}}^K}}
\newcommand{\surfEK}{{F_{\text{surf}}^{K,\sharpEps}}}
\newcommand{\restEK}{{F_\partial^{K,\sharpEps}}}
\newcommand{\fullEK}{{F^{K,\sharpEps}}}

\titleformat*{\section}{\normalsize\bfseries}

\title{Optimization of the branching pattern in coherent phase transitions}
\author{Patrick Dondl, Behrend Heeren, Martin Rumpf} 
\date{\today}

\begin{document}

\maketitle
 
\begin{abstract}
Branching can be observed at the austenite-martensite interface of martensitic phase transformations. 
For a model problem, Kohn and M\"uller~\cite{KoMu92,KoMu94} studied a branching pattern with optimal scaling of the energy with respect to its parameters. 
Here, we present finite element simulations that suggest a topologically different class of branching patterns and derive a novel, low dimensional family of patterns. 
After a geometric optimization within this family, the resulting pattern bears a striking resemblance to our simulation. 
The novel microstructure admits the same scaling exponents but results in a significantly lower upper energy bound. 
\end{abstract}

\begin{center}
 \emph{keywords}: phase transition, finite element simulation, branching pattern, reduced model
\end{center}

\section{Introduction}

In this article we study a non convex scalar variational problem~\cite{Kh83,Ro69,BaJa87,KoMu92} considered to be a model for Austenite-Martensite interfaces in solid-solid phase transitions. Such interfaces often exhibit characteristic fine-scale phase mixtures, 
where the Austenite phase of higher crystallographic symmetry is separated from thin layers of alternating Martensitic phase variants. 
In the seminal work by Kohn and M\"uller~\cite{KoMu92,KoMu94} the energy functional 
\begin{equation}\label{eq:sharp}
\sharpEnergy[u] = \int_\Omega \sharpEps \, |u_{yy}| +  u_x^2 \, \diff x \diff y
\end{equation}
was studied for suitable $u\colon \Omega \to \R$ with $\Omega = (0,L)\times (0,1) \subset \R^2$ 
under the constraint that $u_y = \pm 1$ almost everywhere and for $u=0$ on 
$\partial \Omega$. They show that there are constants $c,\, \upper >0$ such that $c\,\sharpEps^{2/3}L^{1/3} \leq \min \sharpEnergy[u] \leq \upper \, \sharpEps^{2/3}L^{1/3} $. 
In order to recover the upper bound they use an explicitly constructed branching pattern which -- after some optimization -- yields a constant $C = 6.86$. 
In~\cite{Co00} Conti proved self similarity of the energy minimizer near the Dirichlet boundary at $x=0$, 
in the sense that the sequence $u_j(x, y) = \theta^{-2j/3} u(\theta^j x, \theta^{2j/3} y)$ admits $W^{1,2}$-strongly converging subsequence.
A corresponding diffuse interface Landau-type energy functional approximating $\sharpEnergy$
is given by \begin{equation}\label{eq:diffuse}
 \diffEnergy[u] = \int_\Omega \tfrac12 \, \diffEps^2 \, u_{yy}^2 + \penaltyConst \, (u_y^2 -1)^2 +  u_x^2 \, \diff x \diff y\,,
\end{equation}
with given $\sigma>0$, which corresponds to an effective interfacial energy of $\sharpEps = \frac{2\sqrt{2\sigma}}{3} \diffEps$ in~\eqref{eq:sharp}, cf. eq (1.3) in \cite{KoMu94}.
Here, we report on high precision subdivision finite element simulation of this diffuse model which allowed us to observe a particular branching pattern, that differs substantially from the construction studied in~\cite{KoMu92} as well as other proposed patterns~\cite{Li03a}.  
Based on this insight from the numerical experiment, we construct a novel, low dimensional family of branching microstructures for the sharp 
interface model \eqref{eq:sharp}. Optimizing over all geometric degrees of freedom in our model, the resulting pattern shows a striking geometric resemblance of the subdivision finite element simulation results.
Furthermore, our branching pattern yields a significantly lower upper bound constant $\upper$ for the energy.

\section{Optimal branching pattern observed via FEM simulation for diffuse interfaces} \label{sec:FEM}
The branching pattern investigated in this paper is derived from the results of a finite element simulation of the diffuse interface model. 
To minimize the energy \eqref{eq:diffuse} via a conforming finite element method $H^2$--regular ansatz functions have to be taken into account. 
To this end we fix $L=1$ and use a conforming subdivision finite element approach ~\cite{CiLo11, CiOrSc00, DoShBh07} with $C^2$ basis functions which are quartic polynomials 
defined via Loop subdivision on a regular triangular mesh on the domain $\Omega = (0,L)\times(0,1)$. 
We consider zero Dirichlet boundary conditions at $x=0$  and $x=L$ and periodic boundary conditions at $y=0$ and $y=1$.
\begin{wrapfigure}{r}{0.5\columnwidth}
\centering
\includegraphics[width=0.48\columnwidth]{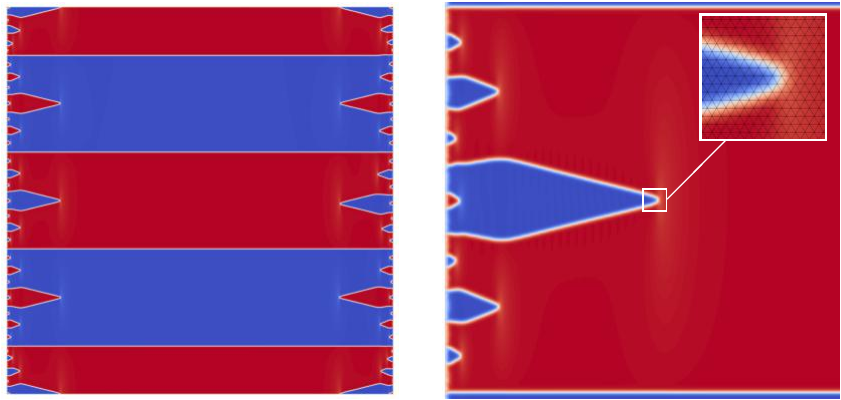}\\[1ex]
\caption{Different zooms are displaced for a discrete minimizer of the energy $\diffEnergy$ in \eqref{eq:diffuse} with 
$10^6$ elements on $\Omega = [0,1]^2$ with $\diffEps = 2\sqrt{5}\cdot 10^{-3}$ and $\sigma = 10$, where blue encodes $u_y \approx -1$ and red $u_y \approx +1$.}
\label{fig:femResult}
\end{wrapfigure}
To compute a discrete minimizer we consider a discrete  gradient flow approach combined with a convex-concave splitting proposed
by D. Eyre~\cite{Ey98} for the Cahn-Hilliard  equation,
which turns out to be unconditionally stable.
By this method the branching pattern can robustly be computed independent of the choosen initial data.
Fig. \ref{fig:femResult} depicts the numerical computed minimizer of the the energy \eqref{eq:diffuse}.
Obviously, the pattern is topologically different from the branching pattern investigated by Kohn and M\"uller \cite{KoMu92}, where the inner needles are not taken into account.
Furthermore, the pattern in our simulations differs geometrically from the pattern proposed by Li \cite{Li03a}, where interior needles are generated on every level whereas in our pattern they arise \emph{every other} level. 
Let us mention that a topologically equivalent, however geometrically somewhat distorted pattern has been found in less accurate simulations by Muite~\cite{Mu09}.

\section{Optimal branching pattern for a reduced sharp interface model}\label{sec:discrete}
Now, we will use the findings for the diffuse interface model from the previous section to derive a geometrically simple, reduced model
for the sharp interface problem~\eqref{eq:sharp},
optimize over its degrees of freedom and compare it in terms of the stored energy with correspondingly optimized version of 
the pattern proposed by  Kohn, M\"uller and by Li.

\begin{figure}[h]
\centering
\includegraphics{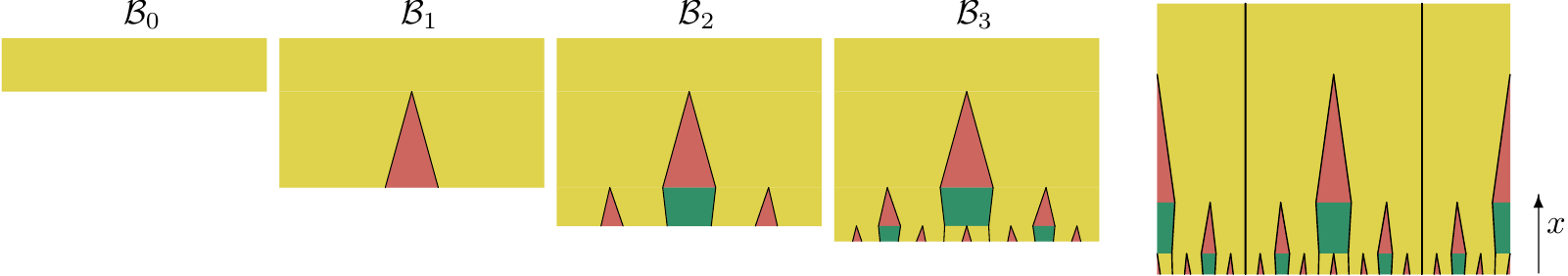}\\[1ex]
\caption{Constructing a periodic cell $\per_3$ with $K\!=\!3$ levels. Starting from the intial word $\branch_0 = \word_0 = \facet$, 
the automaton is applied three times to construct $\branch_1, \branch_2, \branch_3$ (using the (NEW) scheme).
Afterwards, a $y$-periodic cell $\per_3$ is constructed by reflecting half of the pattern $\branch_3$ to either side along two straight additional interfaces (right). 
Color code: facet (yellow), spike (red) and trunk (green), black lines denote interfaces separating regions with $|u_y|= \pm 1$.
All patterns $\branch_k$ and $\per_3$ have been rotated by $90^\circ$ for better visualization. Note: The color coding does \emph{not} represent different phases of martensite (as in Fig. \ref{fig:geometry}). }
\label{fig:topology}
\end{figure}

\begin{figure}[h]
\centering
\includegraphics{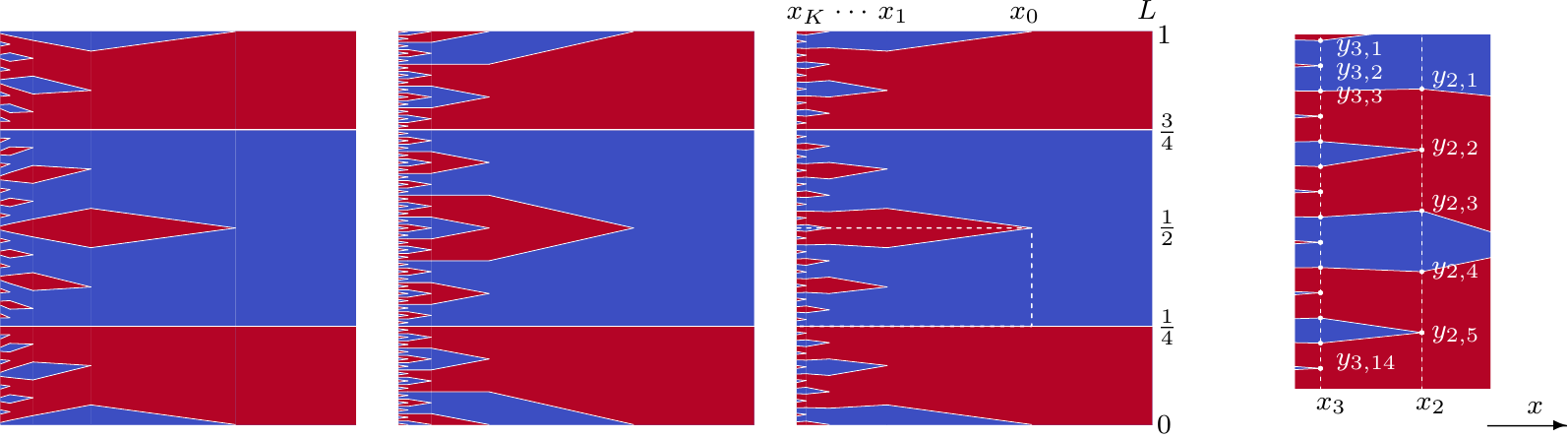}\\[1ex]
\caption{Periodic cell $\per_4$ for (KM), (L) and (NEW) from left to right, with $x_k = \theta^kl$ for $k \geq 0$ and $l < L$. 
Here the color encodes the sign of $u_y$. 
The computational domain $\comp_K = [x_K,x_0] \times [\tfrac14,\tfrac12 ]$ is surrounded by the dotted line in the (NEW) scheme.
By construction, the pattern is symmetric along $y = \tfrac12$, ensuring periodic boundary conditions, i.e. $u(x,0) = u(x,1)$.
Degrees of freedom in $y$-direction are shown on the right (white dots), i.e. $y_{k,n}$ for $n = 1, \ldots, n_k$ for each level $k$.
Later, $\per_K$ is rescaled in $y$-direction by $N^{-1}$ to allow $N$ repetitions in $y$-direction.}
\label{fig:geometry}
\end{figure}

\textbf{Topology.} For the reduced model we restrict the ansatz space for $u$ via an assumption on the geometry of the interfaces where the gradient of $u$ jumps.
In fact, we suppose this interface set to consist of piecewise polygonal lines separating regions with $u_y=\pm 1$.
The branching pattern is now constructed by defining needles of different size, which are bounded by these polygonal lines.
We assume that one needle consists of a trunk (t) and a spike (s) on top of it, material outside a needle is refered to as facet (f).
The generation of a particular branching pattern is described best by the action of an \emph{automaton} $\auto$.
We consider the alphabet $\{\facet,\, \spike,\, \trunk, | \}$, where $|$ denotes an interface.
The initial word is given by $\word_0 = \, \facet\,$. 
A branching pattern after $K$ iterations is defined by  $\branch_K = \word_0; \, \auto(\word_0); \, \ldots; \,\auto^K(\word_0)$, where $;$ represents the generation of a new level.
The  action of $\auto$ on a word is realized letter by letter where $\auto(|) = |$ for all patterns. 
Different branching patterns are characterized by the action of $\auto$ on the three letters $\{ \facet,\, \spike,\, \trunk \}$.
The new branching pattern (NEW) described in Section \ref{sec:FEM}  is characterized by $\auto(\spike) = \trunk$ and $\auto( \trunk) = \auto(\facet) = \facet\, | \,\spike\,| \, \facet \,$ (cf. Fig.\ref{fig:topology}).
The pattern investigated by Kohn and M\"{u}ller~\cite{KoMu92} (KM) is generated by the rules
$\auto(\spike) =  \trunk$, $\auto(\trunk) =  \trunk$ and $ \auto(\facet)=  \facet\,| \, \spike\,| \, \facet$ and the pattern proposed by Li~\cite{Li03a} (L) is given by the rules
$\auto(\spike) =  \auto(\facet) =  \facet\, | \, \spike\,| \, \facet$. 
There are no trunks $\trunk$ in the Li model.
A periodic cell $\per_K$ with $K$ levels is now defined by computing $\branch_K$ (as described above) and then reflecting half of the pattern to either side along two straight additional interfaces (as shown in Fig. \ref{fig:topology}, right).
The constructed pattern consists of connected components (bounded by polygonal lines) which will later represent different phases of martensite, i.e. regions with $u_y=\pm 1$ (cf. Fig. \ref{fig:topology} vs. \ref{fig:geometry}).

\textbf{Geometry.} Now we describe the geometric arrangement of $\Omega = [0,L]\times[0,1]$ (cf. Fig. \ref{fig:geometry}). 
There is a region without branching between $x = l$ and $x =L$.
The $k$th branching level (whose topology is described by the automaton) is located between $x_k = \theta^kl$ and $x_{k-1} = \theta^{k-1}l$, for $k = 1, \ldots , K$.
Finally there is a boundary layer (again without branching) between $x = 0$ and $x_K$.
By construction, the pattern is symmetric along $y = \tfrac12$, ensuring periodic boundary conditions, i.e. $u(x,0) = u(x,1)$.
We choose constant zero boundary condidions $u(x,0) = u(x,1) = 0$ which directly implies $u(x,\tfrac12) = 0$, i.e. $u_x(x,\tfrac12) = 0$.
We assume that the tip of an (initial) needle and the tip of all needles being generated \emph{within} this initial needle lie on a horizontal line. 
Note that this only applies to (NEW) and (L) as (KM) never produces new needles \emph{within} one existing needle.

\textbf{Elastic energy.} The computational domain $\comp_K \subset \per_K$ is limited to $\comp_K = [\theta^Kl,l] \times [\tfrac14,\tfrac12 ]$, as indicated by the dotted lines in Fig. \ref{fig:geometry}.
The geometry of $\comp_K$ is described by $\theta$ and the $y$ coordinates of the interior vertices of the polygonal interfaces.  
We denote the vertices on the line with $x_k=\theta^k l$ from top to bottom by $y_{k,n}$ for $n=1,\ldots n_k$ (cf. Fig. \ref{fig:geometry}, right) 
and the open sets of constant slope $u_y$ between the lines $x_{k-1}=\theta^{k-1} l$ and $x_k=\theta^{k} l$ by $\area_{k,i}$ for $i=1,\ldots i_k$. 
As $u$ is a piecewise linear function by construction, $u_x$ is constant on the regions $\area_{k,i}\,$ for all $i,k$. 
Since $u_x(x,\tfrac12) = 0$, we have $u_x|_{\area_{k,0}} = 0$.
For $[y_{k-1, \overline{n}(k,i)},\, y_{k,\underline{n}(k,i)}]$ being the upper edge of $\area_{k,i}$, the piecewise constant quantity $u_x$ on the stripe $ \comp_K \cap \{ x_{k} \leq \x \leq x_{k-1} \}$ is described iteratively by
\begin{equation}\label{eq:recursiveFormUx}
 u_x|_{\area_{k,i+1}} = u_x|_{\area_{k,i}} + 2 \,\text{sign}\left(u_y\big|_{\area_{k,i}}\right)  \frac{ \y_{k-1, \overline{n}(k,i)} - \y_{k,\underline{n}(k,i)} }{\theta^{k-1}l - \theta^{k}l}\, , \quad i \geq 0\, .
\end{equation}

The elastic energy $\int_{\comp_K} u_x^2 \dxdy$ is given by $\elastE[\comp_K] = \sum_{k=1}^K \sum_{i > 0} ( u_x|_{\area_{k,i}})^2 \cdot | \area_{k,i} |$.
Due to the symmetry of $\per_K$ we have $\elastE[\per_K] = 4\, \elastE[\comp_K]$.
Note that there is no contribution to the elastic energy in the stripe $\{ (x,y) : l < x < L\, \}$.
If we now rescale $\per_K$ in $y$-direction by $N^{-1}$ (to allow $N$ repetitions of $\per_K$) we have $u_x \sim N^{-1}$ and $|\area_{k,i}| \sim N^{-1}$.
The total elastic energy $\elastEK$ is now given by summing over all $N$ repetitions of $\per_K$, i.e. if we collect all $\y_{k,i}$ in a vector $\Y$ we have
\begin{equation}\label{eq:elastEK}
 \elastEK[\theta, N, l, \Y] = \int_{x > x_K} u_x^2 \dxdy = 4N^{-2}\,\sum_{k=1}^K \sum_{i > 0} \Big( u_x|_{\area_{k,i}}\Big)^2 \cdot | \area_{k,i} | \quad \sim \quad c_\text{elast}(\theta,\Y) \,\, l^{-1} N^{-2} \, .
\end{equation}

\textbf{Surface energy.} Let $\numIk$ denote the number of interfaces in $\stripe_k = \{ (x,y) : x_{k} < x < x_{k-1}\} \subset \per_K$.
Then we obtain $\int_{\stripe_k}\sharpEps |u_{yy}| \dxdy = 2\sharpEps \, (x_{k-1}-x_k) \, \numIk = 2\sharpEps \, (1-\theta)\theta^{k-1}l \, \numIk$ for the interface energy which does not depend on the inner vertices $y_{k,i}$. 
Obviously $\numIk = 2 \cdot (1 + 2\sum_{i=1}^k s_i)$, where $s_i$ denotes the number of spikes in $\stripe_i$ between two parallel interfaces shown as vertical black lines in Fig. \ref{fig:topology} (right).
We deduce $s_i = 2^{i-1}$ for (KM) and $s_i = 3^{i-1}$ for (L), whereas for (NEW) the recursive relation $s_1 = 1$, $s_2=2$ and $s_{i} = 2s_{i-1} + s_{i-2}$ for $i>2$ implies $s_i = (\alpha - \beta)^{-1}(\alpha^{i} - \beta^{i})$, where $\alpha = 1 + \sqrt{2}$ and $\beta = 1 - \sqrt{2}$.
Plugging this into the formula for $\numIk$ yields $\numIk = \alpha^{k+1} + \beta^{k+1}$ for (NEW), $\numIk = 2(2^{k+1}-1)$ for (KM) and $\numIk = 2\cdot3^k$  for (L). 
If we take into account the $2$ interfaces in $\{ (x,y) \in \per_K : l < x < L\}$ and consider $N$ repetitions of a rescaled cell $\per_K$ we get 
\begin{equation}\label{eq:surfEK}
  \surfEK[\theta, N, l] = \int_{x > x_K}\sharpEps |u_{yy}| \dxdy = 2 \sharpEps N(1-\theta) l \sum_{k=1}^K \theta^{k-1}\numIk\,  + 4 \sharpEps N (L-l) \quad \sim \quad c_\text{surf}(\theta) \,\, \sharpEps\,l\, N \, .
\end{equation}

\textbf{Boundary conditions.} So far we have neglected the energy in the domain $\Omega_\partial^K = \{ (x,y) : 0 \leq x < x_K \}$.
In particular we have not specified how to fulfill the boundary conditions $u = 0$ at $x=0$ in the reduced model.
As there is no closed form formula for the elastic energy on a stripe $\stripe_k$, we can evaluate (\ref{eq:elastEK}) only for finite $K$.
Similar to the proof of Lemma 2.3 in \cite{Co00} we construct $\tilde u: \Omega_\partial^K \rightarrow \R$ explicitly, with $|\tilde u_y| = 1$, $\tilde u (0,y)  = 0$ and $\tilde u(x_K,y) =u(x_k, y)$ and estimate the elastic energy in $\Omega_\partial^K$ by $\int_{\Omega_\partial^K} {\tilde u}_x^2 \dxdy$.
Furthermore, we will ensure $\int_{\Omega_\partial^K} \sharpEps  |{\tilde u}_{yy}| \dxdy \leq \int_{\Omega_\partial^K} \sharpEps |u_{yy}| \dxdy = 2\sharpEps \, N(1-\theta) l \sum_{k>K} \theta^{k-1}\numIk $, cf. (\ref{eq:surfEK}).
The sum on the right hand side converges if $\alpha\theta < 1$ for (NEW), $2\theta < 1$ for (KM) and $3\theta< 1$  for (L) and the limit can be computed explicitly.
Hence we get a remaining energy
\begin{equation}\label{eq:restEK}
   \int_{\Omega_\partial^K}\sharpEps |u_{yy}| + u_x^2 \dxdy \leq \frac{c_{\text{rem}}(\theta)}{ \theta^Kl N^2}  + \frac{\| u(x_K,.)\|^2_{L^2(0,1)}}{ \theta^Kl } +  2\sharpEps \, N(1-\theta) l \sum_{k>K} \theta^{k-1}\numIk  =: \restEK[\theta, N, l, (\y_{K,i})_i ].
\end{equation}
Note that $\| u(x_K,.)\|^2$ depends on $(\y_{K,i})_i$ and scales like $N^{-2}$.
Furthermore the constant $c_{\text{rem}}(\theta)$ is the limit of a geometric series which converges if $\alpha^2\theta > 1$ for (NEW), $4\theta > 1$ for (KM) and $9\theta > 1$ for (L). \\

\textbf{Optimization.} Using (\ref{eq:elastEK}), (\ref{eq:surfEK}) and (\ref{eq:restEK}) we obtain the objective functional for $\sharpEps > 0$ and $K \in \N$
\begin{equation}\label{eq:fullEK}
   \fullEK[\theta, N, l,  \Y ] = \elastEK[\theta, N, l, \Y] + \surfEK[\theta, N, l]  + \restEK[\theta, N, l, (\y_{K,i})_i ]\, .
\end{equation}
Optimizing $\fullEK$ for $N$ yields $N \sim \sharpEps^{-1/3}\, l ^{-2/3}$, i.e. we get the expected scaling $\fullEK \sim \sharpEps^{2/3}\, l ^{1/3}$.

We now optimize $ \fullEK$ in (\ref{eq:fullEK}) with respect to $\theta$, $l$, $\Y$ for fixed $\sharpEps>0$, $L = 0.5$ and $K \in \N$ for the three different reduced models (L), (KM) and (NEW).
One may also consider $N \in \N$ as an additional degree of freedom which is realized by assuming $N \in \R$ in the simulations.
In oder to compare the results to the finite element simulation of the diffuse interface model shown in Fig. \ref{fig:femResult}, 
where $\diffEps = 2\sqrt{5}\cdot 10^{-3}$ and $\sigma = 10$ we set $\sharpEps = \frac{2\sqrt{2\sigma}}{3}\diffEps \approx 0.013$ in the reduced model; cf. eq (1.3) in \cite{KoMu94}.\\
The simulation is initialized on level $K=4$ with the patterns shown in Fig. \ref{fig:geometry}.
After optimizing all degrees of freedom on this level, $K$ is increased by one by extrapolating the solution of the last level and the whole pattern is optimized again.
This procedure is repeated iteratively.

\section{Discussion and comparison of the branching pattern}\label{sec:result}

\begin{wrapfigure}{r}{0.45\columnwidth}
\centering
\includegraphics{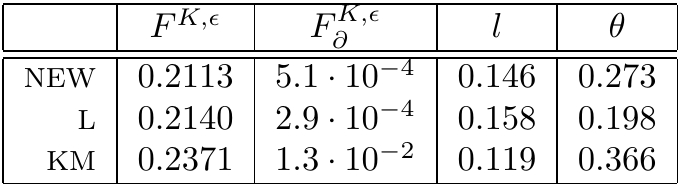}\\[1ex]
\caption{Optimal values for NEW ($K = 17$), L ($K=14$) and KM ($K=12$) for fixed $N=2$.}
\label{fig:optValuesFixed2N}
\end{wrapfigure}

\begin{figure}[h!]
 \centering
\includegraphics{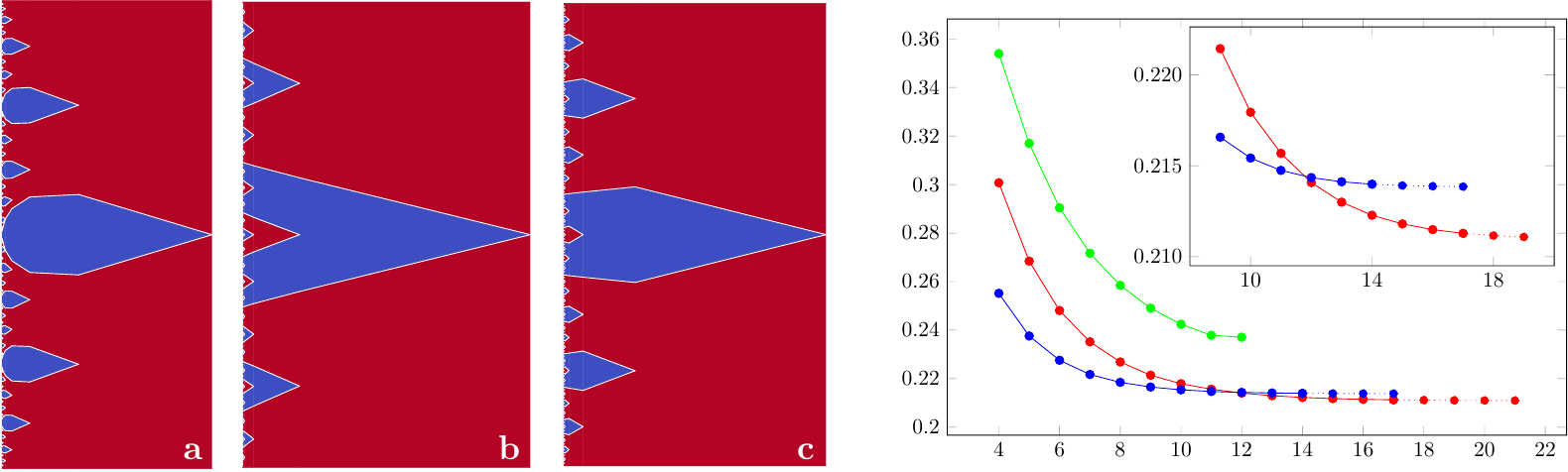}\\[1ex]
\caption{From left to right: optimal patterns of the reduced model (with $\sharpEps = 0.013$ and fixed $N=2$) for (KM) with $K=12$, (L) with $K=14$ and (NEW) with $K=17$, restricted to the domain $[\theta^Kl,l] \times [3/8,5/8]$.
Right: Optimal energy for (NEW) in red, (L) in blue and (KM) in green, with enlargements of crucial details.
Dotted lines are extrapolations without optimization, nevertheless, energies are monotonically decreasing.
Numerical optimization breaks down for (KM) if $K>12$ due to degeneration of geometry on high levels.}
\label{fig:discreteResults}
\end{figure}

The optimal patterns of the reduced models (for fixed $N\!=\!2$) are compared to each other in Fig. \ref{fig:discreteResults}. 
First, the reduced (KM) pattern degenerates already at level $K=13$ , i.e. monotonicity of the sequence $(y_{k,i})_i$ for $k \approx K$ is violated.
In particular, it is energetically not competitive (cf. Fig. \ref{fig:optValuesFixed2N} and \ref{fig:discreteResults}, right).
Qualitatively, the optimal pattern of (NEW) in Fig. \ref{fig:discreteResults}c,  is very similar to the optimal pattern of the diffuse interface model in Fig. \ref{fig:femResult}.
Striking is the fact that all needles of (NEW) are actually diamond-shaped, i.e. bounded by a polygon consisting of \emph{four straight} lines.
Up to some numerical perturbations (probably due to boundary effects) this can also be observed in the pattern in Fig. \ref{fig:femResult}.
However, these straight lines are not mandatory in the reduced model due to the freedom to move $y_{k,i}$. 
For instance, this is not the case for the optimal pattern of (L) and (KM), cf. Fig. \ref{fig:discreteResults}a and \ref{fig:discreteResults}b.
The similarity between the optimal (NEW) pattern and the optimal pattern of the diffuse interface model also holds quantitatively. 
First, for (NEW) the optimal value $\theta = 0.273$ reflects very well the corresponding measured ratio $\theta \approx 0.26$ in Fig. \ref{fig:femResult}.
Second, the ratio "width-to-height" of the largest needle is approx. $2.6$ in Fig. \ref{fig:femResult} and $2.7$ for the optimal pattern of (NEW), Fig. \ref{fig:discreteResults}c.
Furthermore, this ratio decreases in both patterns almost equally when going to the second largest needle.

\begin{wrapfigure}{r}{0.5\columnwidth}
\centering
\includegraphics{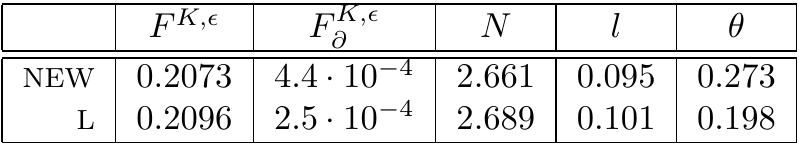}\\[1ex]
\caption{Optimal values for NEW ($K = 17$) and L ($K=14$) with free $N\in\R$. Note that $ \Delta\fullEK \!\sim\! 10^{-3}$ whereas $\restEK\!\sim\! 10^{-4}$. }
\label{fig:optValues}
\end{wrapfigure}

Among the reduced models the (NEW) pattern performs best energetically, with a relative difference of $>1\%$ w.r.t. the optimal energy of (L), see Fig. \ref{fig:optValuesFixed2N}.
Although the difference is not particularly large, it is stable with respect to a variation of the model parameters $\sharpEps$, $L$ and $N$ 
(not shown here).
In particular, this remains true when optimizing over all degrees of freedom (cf. Fig. \ref{fig:optValues}).
Table \ref{fig:optValues} further reveals that the difference $ \Delta\fullEK$ of the optimal energies is of order $10^{-3}$ whereas the error (due to the estimate (\ref{eq:restEK}) in the boundary layer $\Omega_\partial$) is of order $10^{-4}$.
Furthermore, the estimate (\ref{eq:restEK}) is not sharp and hence $\restEK$ is dominant for small $K$ which leads to a substantial over-prediction of the energy of (NEW) for $K < 13$ (cf. Fig. \ref{fig:discreteResults}, right).
To compare the results to the upper constant $\upper = 6.86$ given in \cite{KoMu92} we compute $\upper = \min( \fullEK) L^{-1/3} \, \sharpEps^{-2/3}$ 
and get $\upper = 4.81$ ($N=2$ fixed) or even $\upper = 4.72$ (free $N \in \R$) for (NEW) and $\upper = 4.87$ or $\upper = 4.78$ for (L), respectively.\\

\textbf{Acknowledgements:} B. Heeren and M. Rumpf acknowledge support of the Collaborative Research Centre 1060 and
the Hausdorff Center for Mathematics, both funded by the German Science foundation.  
The authors thank S. Conti, B. Zwicknagl and J. Diermeier for various discussions and helpful comments.\\

\bibliographystyle{plain}
\bibliography{references}

\end{document}